\documentclass[final]{siamltex}
\usepackage{amsfonts}
\usepackage{amssymb}

\newcommand{\xdownarrow}[1]{%
  {\left\downarrow\vbox to #1{}\right.\kern-\nulldelimiterspace}
}
\usepackage{titlesec}
\titleformat{\section}
  {\normalfont\Large\bfseries\centering}{\thesection}{1em}{}

\title{ The normal bundle of a rational curve on a generic quintic threefold}

\author{ B. Wang\\
May 25, 2015}

\begin{document}

\maketitle

\begin{abstract}  This is another proof of the same result in [9].
Let $X_0$ be a generic quintic hypersurface in $\mathbf P^4$ over $\mathbb C$ and
 $c_0$ a regular map $\mathbf P^1\to X_0$ that is generically one-to-one to its image.    In this paper, we show\par
(1) $c_0$ must be an immersion, i.e. the differential $(c_0)_\ast: T_{t}\mathbf P^1 \to T_{c_0(t)} X_0$ is injective\par\quad\hspace{3pt} 
at each $t\in \mathbf P^1$,\par
(2) the normal bundle of $c_0$ satisfies
$$H^1(N_{c_0/X_0})=0.$$

\end{abstract}

\section{Introduction}

Throughout the paper, we work over $\mathbb C$.
Let $X_0$ be a generic quintic threefold in $\mathbf P^4$ over $\mathbb C$. 
Let $c_0: \mathbf P^1\to X_0$ be a  birational map onto its image. 
  The regular map $c_0: \mathbf P^1\to \mathbf X_0$ induces a differential map
\begin{equation}\begin{array}{ccc}
(c_0)_\ast|_t: T_{t}\mathbf P^1 &\rightarrow & T_{c_0(t)}X_0\end{array}
\end{equation}
point-wisely, which induces an injective  morphism on the sheaf module, denoted by $(c_0)_\ast$
\begin{equation}\begin{array}{ccc}
(c_0)_\ast: T_{\mathbf P^1 } &\rightarrow & c_0^\ast(T_{X_0}).
\end{array}\end{equation}

\bigskip

\begin{theorem}  With above set-up, for a generic $X_0$,

(1) $c_0$ is an immersion , i.e. there exists a bundle, called normal bundle, 
$$N_{c_0/X_0}$$ over $\mathbf P^1$ uniquely determined by $c_0$ such that the sequence 
$$\begin{array}{ccccccccc}
0 &\rightarrow & T_{\mathbf P^1} &\stackrel{(c_0)_\ast}\rightarrow &
c_0^\ast(T_{X_0}) &\rightarrow & N_{c_0/X_0}  &\rightarrow & 0,
\end{array}$$
is exact, 
\par
(2) the normal bundle satisfies
\begin{equation} H^1(N_{c_0/X_0})=0. \end{equation}
\par

\end{theorem}

\bigskip

\begin{corollary}
Let $X_0\subset \mathbf P^4$ be a generic quintic hypersurface, and $c_0$ as in theorem 1.1. 
Then the bundle 
$N_{c_0/X_0}$  as in theorem 1.1  has an isomorphism
$$N_{c_0/X_0}\simeq \mathcal O_{\mathbf P^1}(-1)\oplus \mathcal O_{\mathbf P^1}(-1).$$

\end{corollary}
\bigskip

\begin{proof} of corollary 1.2 following from theorem 1.1:  Notice 
$$deg(N_{c_0/X_0})=deg(c_0^\ast(T_{X_0}))-deg(T_{\mathbf P^1})=-2.$$
 It is well-known that the bundle can be split into,
\begin{equation} 
N_{c_0/X_0}\simeq \mathcal O_{\mathbf P^1}(k)\oplus \mathcal O_{\mathbf P^1}(-k-2)
\end{equation}
where $k\geq -1$ is an integer.
By Serre duality
\begin{equation}\begin{array}{cc}
H^1(N_{c_0/X_0})&\simeq H^0( ( N_{c_0/X_0})^\ast\otimes \omega_{\mathbf P^1})\\
& \simeq H^0(\mathcal O_{\mathbf P^1}(-2-k)\oplus \mathcal O_{\mathbf P^1}(k)).\end{array}\end{equation}
By theorem 1.1, $H^1(N_{c_0/X_0})=0$. Hence $-1\leq k\leq -1$. 
Therefore $k=-1$. 

\end{proof}

\bigskip

\subsection{Outline of the proof}\quad\smallskip

The cohomological statement of theorem 1.1 is equivalent to a property of the incidence scheme 
\begin{equation} \Gamma_{X_0}=\{ birational\ to \ its \ image \ maps\ c: \mathbf P^1\to X_0\} \end{equation}  of rational maps  to rational curves of a fixed
degree on generic quintic threefold $X_0$---(1) $\Gamma_{X_0}$ is reduced, (2) it has the expected dimension.  The set of defining equations of this scheme are pretty easy to obtain (see [6]). This property, which is determined by the Jacobian matrix of this set of defining equations   therefore is another expression of theorem 1.1.  Clemens proved that there are components  of $\Gamma_{X_0}$ at whose generic points the Jacobian matrix has full rank ([3], section 1).  But the method can't be used on all components.   In this paper, we prove it for all components. Our general idea of using Jacobian matrices is similar to Clemens', but the detailed steps and the technique are different.   We \par
(I)    replace  the  single  quintic  $X_0$   by  a  generic two parameter 
 family  $\mathbb L$   \par \quad\hspace{3pt} of quintics.\par
(II)  then   show  that the  Jacobian matrix for  the projection $P(\Gamma_{\mathbb L})$   at  a   \par \quad\hspace{3pt}  generic  point has full rank. 
\smallskip

Therefore $P(\Gamma_{\mathbb L})$  is reduced with the expected dimension. Then it follows that the incidence scheme $\Gamma_{X_0}$ is also reduced with
the expected dimension (see [9] for the details). 

\bigskip

By switching $X_0$ to $\mathbb L$, we obtain two free parameters  for the incidence scheme $P(\Gamma_{\mathbb L})$ that come from the deformation of the quintic $f_0$, while the original component $\Gamma_{X_0}$ has no free moduli parameters.
 The manipulation of two free parameters allows us to penetrate the Jacobian matrix. 
The following is the detailed sketch of the proof. For the parameter space of rational maps we use the linear model of moduli maps (used in (1.6) ).  In particular, we do not use a moduli space  of rational maps. By the linear model we mean the affine space $M$,  
\begin{equation} M=\mathbb C^{5d+5}=( H^0(\mathcal O_{\mathbf P^1}(d))^{\oplus 5}\end{equation}
whose open subset parametrizes  the set of non-constant regular maps $$\mathbf P^1\to \mathbf P^4$$ whose push-forward cycles have degree $d$.
 Let $M_d$ be the subset that consists of all generically one-to-one (to its image) maps $c$ whose images $c_\ast(\mathbf P^1)$ have degree $d$. 
Let $S=\mathbf P(H^0(\mathcal O_{\mathbf P^4}(5)))$ be the space of all quintics. 
Let $\mathbb L\subset S$ be an open set of the plane spanned by quintics $f_0, f_1, f_2$.
Let  $$\Gamma_{\mathbb L}\ni (c_0, [f_0])$$ be an irreducible component of the incidence scheme
\begin{equation}\{(c, [f])\subset M_d\times \mathbb L: c^\ast(f)=0\}
\end{equation} 
that is onto $\mathbb L$, where $[f_0]$ denotes the image of $f_0$ under the map
$$\begin{array}{ccc}
H^0(\mathcal O_{\mathbf P^4}(5))-\{0\} &\rightarrow & S. \end{array}$$
We assume $\Gamma_{\mathbb L}$ exists. Let $P$ be the projection $$\Gamma_{\mathbb L}\to M.$$
The idea of the proof is to show that the projection,
$$P(\Gamma_{\mathbb L})\subset M$$ is
a reduced, irreducible quasi-affine scheme of dimension $6$.  The method is straightforward to show its defining polynomials at a generic point have
non-degenerate Jacobian matrix (by that we mean it has full rank). See
definition 1.8 below for a precise definition of a Jacobian matrix.   All differentials and partial derivatives used throughout the paper are in algebraic sense, i.e.
defined as in [7] (because all functions are holomorphic).  In the following we describe its defining polynomials and a differential form representing
the Jacobian  matrix. 
Choose generic $5d+1$ distinct points  $t_i\in \mathbf P^1$ (generic in $Sym^{5d+1}(\mathbf P^1)$). Throughout the paper, unless specified otherwise, we'll use $t_i$ to denote a  complex number  which is a point in an affine open set  $\mathbb C\subset \mathbf P^1$. 
Next we consider differential 1-forms $\phi_i$ on $M$:
  
\begin{equation} \phi_i= d
\left|  \begin{array}{ccc} f_2(c(t_i)) & f_1(c(t_i)) & f_0(c(t_i))\\
f_2(c(t_1)) & f_1(c(t_1)) & f_0(c(t_1))\\
f_2(c(t_2)) & f_1(c(t_2)) & f_0(c(t_2))
\end{array}\right|\end{equation}
for $i=3, \cdots, 5d+1$, and variable $c\in M$, where $|\cdot |$ denotes the determinant of a 
matrix.  Notice $\phi_i$ are uniquely defined provided  the quintics $f_i$ are in an affine open set of $S$, and $t_i\in \mathbb C$ as chosen.
Let \begin{equation} \omega(\mathbb L, \mathbf t)=\wedge_{i=3}^{5d+1}\phi_i \in H^0(\Omega^{5d-1}_M) \end{equation}
be the $5d-1$-form. This $\omega(\mathbb L, \mathbf t)$  is just a collection of all maximal minors  of the Jacobian matrix of  defining polynomials
\begin{equation} 
\left|  \begin{array}{ccc} 
f_2(c(t_i)) & f_1(c(t_i)) & f_0(c(t_i))\\
f_2(c(t_1)) & f_1(c(t_1)) & f_0(c(t_1))\\
f_2(c(t_2)) & f_1(c(t_2)) & f_0(c(t_2))
\end{array}\right| \end{equation}
for the scheme
 $P(\Gamma_{\mathbb L})$, where $\mathbf t=(t_1, \cdots, t_{5d+1})$.

The following proposition asserts the non-degeneracy of the Jacobian matrix of the defining equations of 
$P(\Gamma_{\mathbb L})$. 
  
\bigskip

\begin{proposition} 
The form  $\omega(\mathbb L, \mathbf t)$ is not identically zero on $P(\Gamma_{\mathbb L})$.

\end{proposition}

\bigskip

Then non-degeneracy of the Jacobian matrix means 
\bigskip

\begin{proposition}
If $\omega(\mathbb L, \mathbf t)$ is non-zero  on $P(\Gamma_{\mathbb L})$, the Zariski tangent space of 
$P(\Gamma_{\mathbb L})$ at a generic maximal point must be 
\begin{equation}
dim(M)-deg(\omega(\mathbb L, \mathbf t))
\end{equation}
\end{proposition}

\bigskip

The cohomological statement in theorem 1.1 follows immediately from the propositions on the incidence scheme above.  See [9] for this step.

\bigskip

\begin{proposition}
If propositions 1.3, 1.4 are true, theorem 1.1 is true.
\end{proposition}

\bigskip

Proposition (1.3) is the central part of the proof.  It is a consequence of the study of a Jacobian matrix $\mathcal A(C_M, f_0, f_1, f_2, \mathbf t)$ of a large size $(5d+5)\times (5d+5)$, where $C_M$ stands for local coordinates' system of the space $M$.
In [9], we used the successive blow-ups to study the matrix $\mathcal A(C_M, f_0, f_1, f_2, \mathbf t)$ around a very degenerate point on $P(\Gamma_{\mathbb L})$. In this paper, we avoid the successive blow-ups by directly  studying a generic point of 
$P(\Gamma_{\mathbb L})$.\footnote {But both methods rely on an algebro-geometric process,  ``specialization"}. 

This can be done through a trick.  Let us refer it as a ``break-up trick". This is the process of a sequence of specializations. Roughly speaking, 
we compound the process of breaking up a whole matrix to block matrices, then manipulate the set, $C_M, f_i, \mathbf t$ and the base point $c_g\in P(\Gamma_{\mathbb L})$ to have computable block matrices.  The trick is that we also need to break $$C_M, f_i, \mathbf t, c_g $$
to study each block and
there is no unified $$C_M, f_i, \mathbf t, c_g $$
(generic in some sense) for all block matrices.
But in the end all broken pieces with special sets of $C_M, f_i, \mathbf t, c_g$ must be chosen to coincide at 
the same  generic $$C_M, f_0, f_1, f_2, \mathbf t, c_g.$$
 So specializations must {\bf NOT} be  applied to the entire matrix $\mathcal A(C_M, f_0, f_1, f_2, \mathbf t)$, but they are applied to some  block matrices separately.\par
 
 Let's give a detailed description of  it in the following. 
It suffices to prove the proposition 1.3 for a specific $\mathbb L$. Thus we choose
\begin{equation}
f_0=generic, f_1=z_2z_3z_4 q, f_2=z_0z_1z_2z_3z_4
\end{equation}
where $z_0, \cdots, z_4$ are homogeneous coordinates of $\mathbf P^4$, and $q$ is a generic quadratic, homogeneous polynomial in $z_0, \cdots, z_4$. 
First we write down the differential one form $\phi_i, i=3, \cdots, 5d+1$(expand it using Leibniz rule in differential):
\begin{equation}
\phi_i=d (g_i(c))+\sum_{l=0, j=1}^{l=2, j=2} h_{lj}^i df_l(c(t_j)).
\end{equation}
where $g_i$ is a linear combination of $df_0(c(t_i)), df_1(c(t_i)), df_2(c(t_i))$.
Then it suffices to show the polynomials $$g_i(c)-g_i(c_g),  f_l(c(t_j))-f_l(c_g(t_j)), i=3, \cdots, 5d+1, j=1, 2, l=0, 1, 2$$ 
at a generic point $c_g$ of $P(\Gamma_{\mathbb L})$ form a regular system of parameters for the local ring $\mathcal O_{c_g, M}$ of $M$.
This is the same to show the $(5d+5)\times (5d+5)$ Jacobian matrix $\mathcal A(C_M, f_0, f_1, f_2, \mathbf t)$ 
\begin{equation}
\mathcal A(C_M, f_0, f_1, f_2, \mathbf t)={\partial (g_3, \cdots, g_{5d+1}, f_0(c(t_1)), \cdots, f_2(c(t_2)))\over \partial C_M}
\end{equation}
is non-degenerate at $c_g$.  The following is the trick mentioned above.  For a generic $f_0, q$ and a GENERIC $c_g\in U_\mathbb L$, we can choose a special 
$C_M$ and $\mathbf t$ denoted by $C_M', \mathbf t'$ such that, 
\begin{equation} \mathcal A(C_M', f_0, f_1, f_2, \mathbf t')|_{c_g}\stackrel{Row} \sim 
\left(\begin{array}{cc} I & 0\\
0 & Jac(C_M', c_g)
\end{array}\right)
 \end{equation}
where $I$ is the identity matrix of size $(5d-2)\times (5d-2)$ and $Jac(C_M', c_g)$ is a
$7\times 7$ matrix (this is the break-up of the matrix and it is done in section 3, the step of choosing $C_M$)\footnote{ This break-up requires that 
$c_g$ is birational to its image. If $c_g$ is a multiple cover map, this break-up does not hold.}.  This $Jac(C_M', c_g)$ is the most difficult part in $A(C_M, f_0, f_1, f_2, \mathbf t)$. 
Next to penetrate the $7\times 7$ matrix,   $Ja(C_M', c_g)$, we use  the ``break-up trick" to break the set of $C_M, \mathbb L, \mathbf t, c_g$, i.e. we choose a special $c_g^1\in P(\Gamma_{\mathbb L})$ and another coordinates $C_M''$ to show
that $Jac(C_M'', c_g^1)$ is non-degenerate. Thus $Ja(C_M', c_g)$ is non-degenerate. 
The trick is that those special $c_g1, C_M''$ fail the formula (1.15), therefore should be avoided at the first place(there will be a couple of more similar break-ups of $Jac(C_M'', c_g^1)$ later). 
  Therefore $\mathcal A( C_M, f_0, f_1, f_2, \mathbf t)$ is non-degenerate for any $C_M$, and generic $f_0, f_1, f_2, \mathbf t$.

 \bigskip

\bigskip

\subsection{Technical notations}\quad\smallskip

In this section, we collect all technical notations and definitions used in this paper. Some of them may already be defined before.\bigskip

{\bf Notations}:
\par
(1) $S$ denotes the space all quintics, i.e. $S=\mathbf P(H^0(\mathcal O_{\mathbf P^4}(5)))$.\par
Let $[f]$ denote the image of $f$ under the map
$$\begin{array}{ccc}
H^0(\mathcal O_{\mathbf P^4}(5))-\{0\} &\rightarrow & S. \end{array}$$

(2) Let $$M$$ be
$$\mathbb C^{5d+5}\simeq (H^0(\mathcal O_{\mathbf P^1}(d))^{\oplus 5}$$ 
and $M_{d}$ be the subset that parametrizes  all birational-to-its-image maps $$\mathbf P^1\to \mathbf P^4$$ whose push-forward cycles  have degree $d$. 
 \par
(3) Throughout the paper, if $$c:  \mathbf P^1\to \mathbf P^4, $$
is regular,   $c^\ast(\sigma)$   denotes the pull-back section  of
section $\sigma$ of some bundle over $\mathbf P^4$. The vector bundles will not always be specified, but they are apparent in the context.\par
(4) Let $Y$ be a scheme,  $y\in Y$ be a closed point, $Z\subset Y$ be a subscheme (open or closed) and $\mathcal M$ be a quasi-coherent
sheaf of $\mathcal O_{Y}$-module.
Then $\mathcal O_{y, Y}$ denotes the local ring, $\Omega_{Y}$ denotes the sheaf of differentials, 
  $\mathcal M|_{(Z)}$ denotes the inverse sheaf module $i^\ast(\mathcal M)$ where $i: Z\hookrightarrow Y$ is the embedding. We call  
$\mathcal M|_{(Z)}$ the restriction of $\mathcal M$ to $Z$. 
$\mathcal M|_{Z}$ denotes the localization of $\mathcal M$ at $Z$, which is a $\mathcal O_{Z, Y}$ module.  Thus
$$\mathcal M|_{(\{y\})}=\mathcal M|_Z\otimes k(y),$$
where $k(y)$ is the residue field of the maximal point $\{y\}$.  
\par
If $Y$ is quasi-affine scheme, $\mathcal O (Y)$ denotes the ring of regular functions on $Y$.  

\par

(5) If $Y$ is a scheme, $|Y|$ denotes the induced reduced scheme of $Y$.

\bigskip

\begin{definition}
Let $\Gamma$  be an irreducible component of the incidence scheme 
\begin{equation}\{(c, f)\subset M_d\times \mathbf P(H^0(\mathcal O_{\mathbf P^4}(5))): c^\ast(f)=0\}
\end{equation} 
that dominates $S=\mathbf P(H^0(\mathcal O_{\mathbf P^4}(5)))$. \end{definition}
Let $(c_0, [f_0])\in \Gamma$ be a generic point.  Throughout the paper we assume that such a $\Gamma$ exists. 
\bigskip

{\bf Remark}: The existence of such a $\Gamma$ is equivalent to  the assumption of theorem 1.1: 
$X_0$ is generic. 
\bigskip

 \begin{definition} Let $f_1, f_2 \in H^0(\mathcal O_{\mathbf P^4}(5))$ be two quintics different from $f_0$.
Let $\mathbb L$ be an open set of the plane in 
$$ \mathbf P(H^0(\mathcal O_{\mathbf P^4}(5)))$$ spanned by $[f_0], [f_1], [f_2]$ and centered around $[f_0]$.

\bigskip

 Let
\begin{equation} 
\Gamma_{\mathbb L}=\Gamma\cap ( M\times \mathbb L)
\end{equation} be an irreducible component of the restriction of $\Gamma$ to $M\times \mathbb L$ such that it is onto ${\mathbb L}$, and
\begin{equation} 
\Gamma_{f_0}, \ for\ generic\ f_0\in \mathbb L
\end{equation} is an irreducible component of
$$P(\Gamma\cap ( M\times \{[f_0]\}))$$
where $P$ is the projection to $M$.

 \end{definition}

\bigskip

\begin{definition} 
Let $\Delta^n\subset \mathbf C^n$ be an analytic open set with coordinates $x_1, \cdots, x_n$, Let $f_1, \cdots, f_m$ be holomorphic functions on $\Delta$. For any positive integers $m'\leq m, n'\leq n$ and a point $p\in \Delta^n$, 
we define

\begin{equation} \begin{array}{c} 
{\partial (f_1, f_2, \cdots, f_{m'})\over \partial (x_1, x_2, \cdots, x_{n'})}|_{p}=\left (\begin{array}{cccccc}
{\partial f_1\over \partial x_1}  &{\partial f_1\over \partial x_2}   & \cdots & {\partial f_1\over \partial x_{n'}}   \\
{\partial f_2\over \partial x_1}  &{\partial f_2\over \partial x_2} & \cdots & {\partial f_2\over \partial x_{n'}}  \\
 \vdots & \vdots & \cdots &\vdots \\
{\partial f_{m'}\over \partial x_1}  &{\partial f_{m'}\over \partial x_2} & \cdots & {\partial f_{m'}\over \partial x_{n'}}
\end{array}\right)|_p. \end{array}\end{equation}
to be the Jacobian matrix of functions $f_1, \cdots, f_{m'}$ in $x_1, \cdots, x_{n'}$. 
  If $n'=n\geq m'=m$, the differential form $\wedge_{i=1}^m df_i$ is just the collection of all $m\times m$
minors of the Jacobian matrix. 
\par

This Jacobian matrix is just the differential of the composition map
\begin{equation}\begin{array}{c}
(x_1, \cdots, x_{n'}) \\
\downarrow\\
 (x_1, \cdots, x_{n'}, x_{n'+1}, \cdots, x_n)\\
\downarrow\\
 (f_1, \cdots, f_m) \\
\downarrow\\ (f_1, \cdots, f_{m'}),\end{array}\end{equation}
where the first map is defined by the embedding  via the coordinates of $p$.

\end{definition}

\bigskip
This definition depends on all coordinates $x_1, \cdots, x_n$ and  it is crucial.  One of main difficulties of this paper is to search for such coordinates that would make Jacobian matrices simpler.  

\bigskip

In section 2, we prove that original Clemens' conjecture follows from theorem 1.1.  In section 3, we prepare the analytic
coordinates of $M$ for the computation. In section 4, we use
the sheaf of differentials to show the non-vanishing property of $5d$-$1$-form $\omega(M, \mathbf t)$  on the scheme $\mathbf P(\Gamma_{\mathbb L})$. This is the central section of the paper. 
It leads the proof of propositions 1.3, 1.4. Section 5 collects two known examples which emphasize on the  singular rational curves. \bigskip

{\bf Acknowledgment}. We would like to thank Bruno Harris who clears our understanding of the map $(c_0)_\ast$ in (1.2).

\bigskip

\bigskip
\section{Clemens' conjecture}
\quad

 Rational curves on hypersurfaces have been great interests for many years in algebraic geometry. The Clemens' conjecture sits in the center
of many major problems in this area.   In [2], its original 1986 statement, Clemens 
proposed:\par``(1) the generic quintic threefold $V$ admits only finitely many rational curves of each degree.\par
(2) Each rational curve is a smoothly embedded $\mathbf P^1$ with normal bundle
\begin{equation}
\mathcal O_{\mathbf P^1}(-1)\oplus  O_{\mathbf P^1}(-1).\end{equation}\par
(3) All the rational curves on $V$ are mutually disjoint. The number of rational curves of degree $d$ on $V$ is
\begin{equation}
(interesting \ number)\cdot 5^3\cdot d. "
\end{equation}
\smallskip

During the last thirty years, there are many articles on the conjecture. The most of them followed the early idea of Katz ([6]) to show that there is only one irreducible component of the incidence scheme, containing a smooth rational curve and  dominating the space of quintics.      
In 1995, Vainsencher  found the degree 5, 6-nodal rational curves in the generic quintic threefolds ([8]). This partially disproved part (2) in the Clemens' conjecture and leave
the part (1) unanswered. At the meantime Mirror symmetry came to the stage to redefine the approach in part (3). 
Based on Vainsencher's result, in 1999, motivated by the Gromov-Witten invariants in the mirror symmetry, 
Cox and Katz  modified the Clemens' original conjecture to the most current form ([4]):\par
`` Let $V\subset \mathbf P^4$ be a generic quintic threefold. Then for each degree $d\geq 1$, we have\par

(i) There are only finitely many irreducible rational curves $C\subset V$ of degree $d$.\par
(ii) These curves, as we vary over all degree, are disjoint from each other.\par
(iii) If $c: \mathbf P^1\to C$ is the normalization of an irreducible rational curve $C$, then the normal bundle
has isomorphism
$$N_{c/V}\simeq \mathcal O_{\mathbf P^1}(-1)\oplus \mathcal O_{\mathbf P^1}(-1). "$$

\bigskip

{\bf Remark}. Cox and Katz's conjecture (iii) should be understood as in two steps. First  $N_{c/V}$ must be a locally free sheaf, secondly
$$N_{c/V}\simeq \mathcal O_{\mathbf P^1}(-1)\oplus \mathcal O_{\mathbf P^1}(-1). $$
We proved the first by showing that $c_0$ is an immersion. \par
The conjecture is proved to be correct for $d\leq 9$ by the work of Katz ([6]), Johnsen and Kleiman ([5]), and Cox and Katz ([4]),  etc.

\subsection{A proof of Clemens' conjecture}\quad\smallskip

\bigskip
  Clemens' conjecture follows from Theorem 1.1 and corollary 1.2  because the corollary below\bigskip

\begin{corollary}
Let $X_0\subset \mathbf P^4$ be a generic quintic threefold. Then for each degree $d\geq 1$, we have\par
(i) there are only finitely many irreducible rational curves $C_0\subset X_0$ of degree $d$. 
\par
(ii) Each rational curve in (i) is an immersed rational curve with normal bundle
$$N_{c_0/X_0}\simeq \mathcal O_{\mathbf P^1}(-1)\oplus \mathcal O_{\mathbf P^1}(-1).
$$
By  ``immersed rational curve" we mean that the normalization map is an immersion. 
\end{corollary}
\bigskip

\begin{proof} of corollary 2.1 following from theorem 1.1 and corollary 1.2:  The existence of rational curves on a generic quintic was proved in [3], [6].
So it suffices to prove the finiteness. 
Part (i) follows from part (ii). So let's prove part (ii).  Let $C_0$ be an irreducible rational curve of degree $d$ on $X_0$. Then we take a normalization of $C_0$, and denote it by $c_0:\mathbf P^1\to X_0$. Since $X_0$ is generic, we have the set-up for corollary 1.2.  Applying corollary 1.2, we obtain part (ii).    
\end{proof}

 \bigskip

Corollary 2.1  proves the modified Clemens' conjecture, namely parts (i) and (ii)
of Cox and Katz's statements.  Clemens' original conjecture must be modified in the light  of  Vainsencher's result.

\bigskip

\section{Space of rational curves, $M$}

\bigskip
The basis of this paper is the linear model of stable moduli, which begins with the projectivization $\mathbf P(M)$. The space $M$ is an affine space
$\mathbb C^{5d+5}$, therefore is very simple.  But we are interested in some subschemes which are not trivial at all. Our idea is to introduce various
analytic coordinates of each copy $\mathbb C^{d+1}$ in $\mathbb C^{5d+5}$. The purpose of these coordinates is to 
provide various parameters for the local ring $\mathcal O_{c, M}$ so that the Jacobian matrices under these coordinates are either
diagonal or triangular. 
In this section we introduce a couple of coordinates systems $C_M', C_M''$ that will be used for our ``break-up trick".   
\bigskip

However readers may skip this section because without section 4 technical preparation here may seem to be aimless. \bigskip

Let $c_g=(c_g^0, \cdots, c_g^4)\in M_d$ with $$c_g^i\in H^0(\mathcal O_{\mathbf P^1}(d))-\{0\}, i=0, \cdots, 4.$$
We may assume $t\in \mathbb C\subset \mathbf P^1$. 
Because $c_g\in M_d$, we assume $c_g^i(t)=0, i=0, \cdots, 4$ have  $5d$ distinct  zeros $$ \tilde\theta_i^j, for\ i\leq 4.$$

Then each component, $ H^0(\mathcal O_{\mathbf P^1}(d))$ of 
$$ M=H^0(\mathcal O_{\mathbf P^1}(d))^{\oplus 5}$$
 has local analytic ``polar" coordinates \begin{equation}
r_i, \theta_i^j, j=1, \cdots, d, for \ r_i\neq 0\end{equation}
 (for each $i=0, 1, 2, 3, 4$) around
$c_g$ such that
\begin{equation} c^i(t)=r_i \prod _{j=1}^d (t-\theta_i^j).  
\end{equation}
 
Let  coordinates values for $c_g$ be
$$r_l=y_l, \theta_i^j=\tilde\theta_i^j,  l=0, \cdots, 4,  i=0, \cdots, 4,  j=1, \cdots, d.$$
Let $q$ be a generic, homogeneous quadratic polynomial in $z_0, \cdots, z_4$.

Let 
\begin{equation} 
h(c, t)=\delta_1q(c(t))+\delta_2 c_3(t) c_4(t).
\end{equation}
for $c\in M$, where $\delta_i, i=1, 2$ are two none zero complex numbers.
   Let $\beta_1, \cdots, \beta_{2d}$ be the zeros of 
$h(c_g, t)=0$. Also let  
$$h(c, t)=\xi \prod _{i=1}^{2d}(t-\epsilon_i).$$
It is clear that $$\begin{array}{c}
\xi=\delta_1 q(r_0, r_1, r_2, r_3, r_4)+\delta_2 r_3r_4, \ and\\
\epsilon_i\ are \ analytic\ functions\ of\ c.
\end{array}$$
Let the corresponding value of $\xi$ at $c_g$ be $\xi^0$. 
By the genericity of $q$, we may assume $\beta_i, i=1, \cdots, 2d$ are distinct and non-zeros . 
Furthermore we assume $\beta_i$ are distinct for $q=z_1z_2$ and generic $\delta_i$.

\bigskip

\begin{proposition}
Let $U_{c_g}\subset M$ be an analytic neighborhood of $c_g$. 

\par
(a) Let \begin{equation}\begin{array}{ccc}
\varrho: U_{c_g} &\rightarrow &  \mathbb C^{5d+5}
\end{array}\end{equation}
be a regular map that is defined
by 
\begin{equation}\begin{array}{cc} &
 (\theta_0^1, \cdots, \theta_4^d, r_0, r_1, r_2, r_3, r_4) \\
 &\xdownarrow{0.3cm}\scriptstyle{\varrho} \\ &
(\theta_0^1, \cdots, \theta_2^d, \epsilon_1, \cdots, \epsilon_{2d}, 
r_0, \cdots, r_3, \xi).
\end{array}\end{equation} 

Then $\varrho$ is an isomorphism to its image.

\par
(b) Let \begin{equation}\begin{array}{ccc}
\varrho': U_{c_g} &\rightarrow &  \mathbb C^{5d+5}
\end{array}\end{equation}
be a regular map that is defined
by 
\begin{equation}\begin{array}{cc} &
(\theta_0^1, \cdots, \theta_4^d, r_0, r_1, r_2, r_3, r_4) \\
 &\xdownarrow{0.3cm}\scriptstyle{\varrho'}\\ &
(\theta_0^1, \cdots, \theta_2^d, \epsilon_1, \cdots, \epsilon_{2d}, 
r_0, \cdots, r_3, r_4).
\end{array}\end{equation} 
Then $\varrho'$ is an isomorphism to its image.
\end{proposition}
\bigskip

\begin{proof}
It suffices to prove the differential of $\varrho$ at $c_g$ is an isomorphism for a SPECIFIC $q$. So we assume  that 
$$ \delta_1=\delta_2=1, q=z_1z_2$$
This is a straightforward calculation of the Jacobian determinant of $\varrho$. We may still assume that $\beta_i, i=1, \cdots, 2d$ are distinct. Using
the composition of two isomorphisms, we obtain that 
the Jacobian determinant
\begin{equation}
det ({\partial (\tilde \theta_0^1, \cdots, \tilde \theta_2^d, y_0,\cdots, y_3, \xi^0,   \beta_1,\cdots, \beta_{2d})
\over \partial (\theta_0^1, \cdots, \theta_2^d, r_0, r_1, r_2, r_3, r_4, \theta_3^1, \cdots, \theta_4^d)})
\end{equation} is equal to 

\begin{equation}
a\cdot {\partial \xi\over \partial r_4}|_{c_g} \cdot J
\end{equation}
where $a$ is some non-zero number, ${\partial \xi\over \partial r_4}|_{c_g}$ is also non-zero and $J$ is another Jacobian

\begin{equation}
J=\left|\begin{array}{ccc}
{\partial  h(c, \beta_1)\over \partial \theta_3^1}  &\cdots  & {\partial h(c, \beta_{1})\over \partial\theta_4^d }\\
\vdots & \vdots & \vdots\\
{\partial  h(c, \beta_{2d})\over \partial \theta_3^1}  &\cdots  & {\partial h(c, \beta_{2d})\over \partial\theta_4^d }
\end{array}\right|_{\tilde c_2}
\end{equation}
Let $T_i, i=0, d$   be the determinant
\begin{equation}
\left|\begin{array}{ccc}
\beta_{i+1}  & \cdots & \beta_{i+1}^d\\
\vdots & \vdots & \vdots \\
\beta_{i+d}  & \cdots & \beta_{i+d}^d
\end{array}\right|.
\end{equation}
Then we compute the determinant to have 
\begin{equation}
J=(-1)^d T_0T_d \prod_{i=1}^d (c_g^3(\beta_{d+i})c_g^4(\beta_{i}) -c_g^3(\beta_{i}) c_g^4(\beta_{d+i})).\end{equation}
Since $\beta_i$ are distinct and non-zeros, $$T_0\neq 0, T_d\neq 0.$$
Since ${c_g^3(t)\over c_g^4(t)}$ is a non-constant rational function and
$$deg(c_g^3(t))=deg(c_g^4(t))=d$$
then we can always arrange the indexes of $\beta_i$ such that each number 
\begin{equation} ({c_g^3(\beta_{d+i})\over c_g^4(\beta_{d+i})}-{c_g^3(\beta_{i})\over c_g^4(\beta_{i})})
\end{equation} is not zero. Hence 
$$\prod_{i=1}^d (c_g^3(\beta_{d+i})c_g^4(\beta_{i}) -c_g^3(\beta_{i}) c_g^4(\beta_{d+i}))\neq 0.$$

Thus $J$ is non-zero.   Therefore 
\begin{equation}
det({\partial (\tilde \theta_0^1, \cdots, \tilde \theta_2^d, y_0,\cdots, y_3, \xi^0,   \beta_1,\cdots, \beta_{2d})
\over \partial (\theta_0^1, \cdots, \theta_2^d, r_0, r_1, r_2, r_3, r_4, \theta_3^1, \cdots, \theta_4^d)})\neq 0
\end{equation} 
 The proof of part (b) is the same as for part (a). 
We complete the proof. 

\end{proof}

\bigskip

\begin{definition}
By proposition 4.1,  both
\begin{equation} \theta_0^1, \cdots, \theta_2^d, r_0, \cdots, r_3, \xi, \epsilon_1, \cdots, \epsilon_{2d}
\end{equation} and 
\begin{equation} \theta_0^1, \cdots, \theta_2^d, r_0, \cdots, r_4, \epsilon_1, \cdots, \epsilon_{2d}
\end{equation} 

are  local analytic coordinates of $M$ around $c_g$, and
$c_g$ corresponds to the coordinate values
\begin{equation}\begin{array} {c}
 \theta_i^j=\tilde \theta_i^j, i=0, 1,2, j=1, \cdots, d\\
r_l=y_l\neq 0, l=0, \cdots, 4\\
\epsilon_i=\beta_i, i=1, \cdots, 2d
\end{array}\end{equation}
and $\xi^0$. 
\bigskip

For the simplicity, we refer the first coordinates' system \begin{equation} \theta_0^1, \cdots, \theta_2^d, r_0,\cdots, r_3,  \xi, \epsilon_1, \cdots, \epsilon_{2d}
\end{equation}  as $C_M'$ and the second one
\begin{equation} \theta_0^1, \cdots, \theta_2^d, r_0, \cdots, r_4, \epsilon_1, \cdots, \epsilon_{2d}
\end{equation}
as $C_M''$. 
\end{definition}
\bigskip

The following lemma is also a local expression for the calculation later. 
Choose  homogeneous coordinates  $[z_0, \cdots, z_4]$ for $\mathbf P^4$. 
Let 
\begin{equation} 
f_3=z_0z_1z_2(\delta_1q+ \delta_2z_3 z_4).
\end{equation}
where $\delta_i$ are two non-zero complex numbers, and $q$ is a generic, quadratic homogeneous polynomial   in $z_0, \cdots, z_4$. 
Let 
$c_g\in M_d$ as above.
$$f_3(c_g(t))\neq 0.$$
We denote the zeros of $c_g^i(t)=0$ by $\tilde\theta_i^j$ and zeros of
\begin{equation}
(\delta_1q+ \delta_2z_3 z_4|_{c_g(t)})=0
\end{equation}
by $\beta_i, i=1, \cdots, 2d$. 
We assume $\tilde \theta_i^j, i=0, \cdots, 4, j=1, \cdots, d$ are distinct, and $\beta_i, i=1, \cdots, 2d$ are also distinct. 

\bigskip

\begin{lemma}   Recall 
in definition 4.2, 
$$ \theta_0^1, \cdots, \theta_2^d, r_0, \cdots, r_3, \xi,  \epsilon_1, \cdots, \epsilon_{2d}
$$ are analytic coordinates of $M$ around the point $c_g$.

Then \par

(a) the Jacobian matrix

\begin{equation}\begin{array} {c} J( \tilde c_2)\\
\|\\
\left(\begin{array}{cccccc}
{\partial f_3(c_g(\tilde \theta_0^1))\over \partial \theta_0^1}  &\cdots &
{\partial f_3(c_g(\tilde \theta_0^1))\over \partial \theta_2^d} & {\partial f_3(c_g(\tilde \theta_0^1))\over \partial \epsilon_1} &\cdots &
{\partial f_3(c_g(\tilde \theta_0^1))\over \partial \epsilon_{2d}}\\
{\partial f_3(c_g(\tilde \theta_0^2))\over \partial \theta_0^1}  &\cdots &
{\partial f_3(c_g(\tilde \theta_0^2))\over \partial \theta_2^d} & {\partial f_3(c_g(\tilde \theta_0^1))\over \partial \epsilon_1} &\cdots &
{\partial f_3(c_g(\tilde \theta_0^2)\over \partial \epsilon_{2d}}\\
{\partial f_3(c_g(\tilde \theta_0^3))\over \partial \theta_0^1}  &\cdots &
{\partial f_3c_g(\tilde \theta_0^3))\over \partial \theta_2^d} & {\partial f_3(c_g(\tilde \theta_0^3))\over \partial \epsilon_1} &\cdots &
{\partial f_3(c_g(\tilde \theta_0^3))\over \partial \epsilon_{2d}}\\
\vdots&\vdots &\vdots &\vdots &\vdots&\vdots\\
\vdots&\vdots &\vdots &\vdots &\vdots&\vdots\\
{\partial f_3(c_g(\tilde \theta_4^d))\over \partial \theta_0^1}  &\cdots &
{\partial f_3(c_g(\tilde \theta_4^d))\over \partial \theta_2^d} & {\partial f_3(c_g(\tilde \theta_4^d))\over \partial \epsilon_1} &\cdots &
{\partial f_3(c_g(\tilde \theta_0^d))\over \partial \epsilon_{2d}}
\end{array}\right)\end{array}\end{equation}
is equal to 
a diagonal matrix $D$ whose diagonal entries are
\begin{equation}
{\partial f_3(c_g(\tilde \theta_0^1))\over \partial \theta_0^1}, \cdots, {\partial f_3(c_g(\tilde \theta_2^d)\over \partial \theta_2^d}, 
{\partial f_3(c_g(\tilde \theta_3^1))\over \partial \epsilon_1}, \cdots, 
{\partial f_3(c_g(\tilde \theta_4^d))\over \partial \epsilon_{2d}}
\end{equation}
which are all non-zeros.\par
(b) For $i=0, \cdots, 4$, $j=1, \cdots, d$, $l=0, \cdots, 3$
$${\partial f_3(c_g(\tilde \theta_i^j))\over \partial r_l}={\partial f_3(\tilde c_2(t_i))\over \partial \xi}=0.$$

\end{lemma}

\bigskip

\begin{proof}  Note $\tilde \theta_i^j, i=0, \cdots, 4, j=1, \cdots, d$ are distinct and $\beta_i, i=0, \cdots, 2d$ are also distinct.  Thus the coordinates
in definition 4.2 exist. 
 We can rewrite 
\begin{equation}
f_3(c(t))=y \prod_{j=1}^{5d}(t-\alpha_j)
\end{equation}
around $\tilde c_2$.
Then $y, \alpha_j$ are all functions of the analytic coordinates in definition 4.2, 
\begin{equation}
\theta_i^j, \epsilon_1, \cdots, \epsilon_{2d}, y_l.
\end{equation}
More specifically $\theta_i^j, i\leq 2$ are exactly the $3d$ roots, $\alpha_i, i=1, \cdots, 3d$ and $\epsilon_i, i=1, \cdots, 2d$ are 
$\alpha_{3d+i}$, and $y$ is an analytic function of $y_l, l=0, \cdots, 4$. Hence 
\begin{equation} f_3 (c(t))=r_0 r_1r_2 \xi \prod_{i=0, j=1, l=1}^{i=2, j=d, l=2d }(t-\theta_i^j)(t-\epsilon_l).\end{equation}
Both parts of lemma 4.3  follow from the expression (4.26).
We complete the proof. 
\end{proof}

\bigskip

\section{Differential sheaf }
\bigskip

In this section, we prove theorem 1.1, i.e. \bigskip

\begin{equation} H^1(N_{c_0/X_0})=0\end{equation}
at generic $(c_0, [f_0])\in \Gamma$.

\bigskip

\bigskip

\subsection{Non-vanishing $5d-1$-form $\omega(\mathbb L, \mathbf t)$ }
\quad\bigskip

\begin{lemma} 

The $5d$-$1$ form $\omega(\mathbb L, \mathbf t)$ defined in (1.10) is a non-zero form when it is evaluated at generic points of $P(\Gamma_{\Bbb L})$, i.e.
the reduction $\bar \omega(\mathbb L, \mathbf t)$ in the module, 
$$H^0(\Omega_M\otimes \mathcal O_{P(\Gamma_{\Bbb L})})$$ is non zero.

\end{lemma}
\bigskip

This lemma is  proposition 1.3 in the introduction.
\bigskip

It suffices to prove lemma 4.1 for special choices of $f_0, f_1, f_2$ and $t_1, \cdots, t_{5d+1}$.  
So let $z_0, z_1, \cdots, z_4$ be  general homogeneous coordinates of
$\mathbf P^4$. Let $$f_2=z_0z_1z_2z_3z_4.$$  
Let 
$$f_1=z_0z_1z_2 q,$$
where $q$ is a generic quadratic homogeneous polynomial in $z_0, \cdots, z_4$. 
   Choose another generic $f_0$.  
Let $$c_g\in P(\Gamma_{\mathbb L})$$
be a generic in $P(\Gamma_{\mathbb L})$.  By the genericity of $f_0$, 
we may assume 
$c_g=(c_g^0, \cdots, c_g^4)$ such that  $c_g^i\neq 0$ and
$c_g^i(t)=0, i=0, \cdots, 4$ have $5d$ distinct zeros $\tilde\theta_i^j\in \mathbf P^1$. 
 To choose $5d$ points $t_i$ on $\mathbb C\subset \mathbf P^1$, we let \par
(1) $t_1, t_2, t_{5d+1}$ be general among all $(t_1, t_2, t_3)$ satisfying
\begin{equation}
\left|  \begin{array}{cc} f_2(c_g(t_1)) &  f_1(c_g(t_1))\\
f_2(c(t_2)) & f_1(c(t_2))\end{array}\right|=0,
\end{equation}

 \par
(2) $t_3, \cdots, t_{5d}$ be the $5d-2$ complex numbers
\begin{equation}\begin{array}{c}
\tilde \theta_i^j, (i, j)\neq (0, 1), (1, 1), i\leq 2\\
\beta_i, i=1, \cdots, 2d
\end{array}\end{equation}
where $\beta_i$ are the zeros of 
\begin{equation}
 \delta_1 q(c_g(t))+ \delta_2 z_3z_4|_{c(t)}=0.
\end{equation}
where \begin{equation}\begin{array} {c} \delta_1=\left|  \begin{array}{cc} f_0(c_g(t_1)) &  f_2(c_g(t_1))\\
f_0(c_g(t_2)) & f_2(c_g(t_2))\end{array}\right|, \\ \quad \\
\delta_2=\left|  \begin{array}{cc} f_1(c_g(t_1)) &  f_0(c_g(t_1))\\
f_1(c_g(t_2)) & f_0(c_g(t_2))\end{array}\right|.\end{array}\end{equation}

To simply put it, $t_3, \cdots, t_{5d}$ are just the zeros of
\begin{equation}
\delta_1 f_1(c(t))+\delta_2 f_2(c(t))=0.
\end{equation}
excluding two zeros $\tilde\theta_0^1, \tilde\theta_1^1$. 

We claim that 
\begin{equation}  \delta_1\neq 0, \quad 
\delta_2 \neq 0.
\end{equation}
This is because $c_g$ lies in a plane $\mathbb L$, but does not lie in the pencils
$span(f_0, f_1), span(f_0, f_2)$.  Thus 
$$\{(f_0(c_g(t)), f_1(c_g(t)))\}_{t\in \mathbf P^1}$$ span
$\mathbb C^2$.  This implies $$\left|  \begin{array}{cc} f_1(c_g(t_1)) &  f_0(c_g(t_1))\\
f_1(c_g(t_2)) & f_0(c_g(t_2))\end{array}\right|\neq 0$$
Similarly $$\left|  \begin{array}{cc}  f_0(c_g(t_1)) &  f_2(c_g(t_1))\\
f_0(c_g(t_2)) & f_2(c_g(t_2))\end{array}\right|\neq 0.$$

\bigskip

We expand $\phi_i, i=3, \cdots, 5d$ to obtain that
\begin{equation}\begin{array}{cc}  \phi_i=\left|  \begin{array}{cc} f_0(c(t_1)) &  f_2(c(t_1))\\
f_0(c(t_2)) & f_2(c(t_2))\end{array}\right| df_1(c(t_i))+\left|  \begin{array}{cc} f_2(c(t_1)) &  f_1(c(t_1))\\
f_2(c(t_2)) & f_1(c(t_2))\end{array}\right| df_0(c(t_i))& \\  +\left|  \begin{array}{cc} f_1(c(t_1)) &  f_0(c(t_1))\\
f_1(c(t_2)) & f_0(c(t_2))\end{array}\right| df_2(c(t_i))+ \sum_{l=0, j=1}^{l=2, j=2} h_{lj}^i(c_g) df_l(c(t_j))
\end{array}\end{equation}

By the assumption for $t_1, t_2$, $$\left|  \begin{array}{cc} f_2(c_g(t_1)) &  f_1(c_g(t_1))\\
f_2(c_g(t_2)) & f_1(c_g(t_2))\end{array}\right|=0.$$
We obtain  
 \begin{equation}\begin{array}{c}  \phi_i|_{c_g}=\delta_1 df_1(c(t_i))+\delta_2 df_2(c(t_i))+ \sum_{l=0, j=1}^{l=2, j=2} h_{lj}^i(c_g) df_l(c(t_j))\\
=d f_3(c(t_i))+ \sum_{l=0, j=1}^{l=2, j=2} h_{lj}^i(c_g) df_l(c(t_j))
\end{array}\end{equation}
where $$ f_3=\delta_1 f_1+\delta_2 f_2.$$
Notice $\delta_1\neq 0\neq \delta_2$.

\bigskip

To show lemma 4.1, it suffices to show the local holomorphic functions
\begin{equation}\begin{array}{c}
f_3(c(t_3)), \cdots, f_3(c(t_{5d+1})),\\
f_0(c(t_1)), f_1(c(t_1)), f_2(c(t_1)),\\
f_0(c(t_2)), f_1(c(t_2)), f_2(c(t_2)).\end{array}
\end{equation}
are the parameters of the local ring $\mathcal O_{c_g, M}$.

Let \begin{equation}
\mathcal A(C_M, f_0, f_1, f_2, \mathbf t')
\end{equation}
be the Jacobian matrix of functions in (4.10) under an analytic coordinate's system $C_M$ at $c_g$.

Then the lemma 4.1 follows from the following lemma

\bigskip

\begin{lemma}  The  $(5d+5)\times (5d+5)$ matrix \begin{equation}
\mathcal A(C_M, f_0, f_1, f_2, \mathbf t')
\end{equation}
is non-degenerate.

\end{lemma}
\bigskip

\begin{proof} of lemma 4.2:   Let's choose $C_M$ to be $C_M'$ defined in 4.2.
Recall $C_M'$ has coordinates
\begin{equation}\begin{array} {c}
\theta_i^j, i\leq 2\\
\epsilon_i, i=1, \cdots, 2d\\
r_0, \cdots, r_3, \xi
\end{array}\end{equation}

The matrix $\mathcal A(C_M, f_0, f_1, f_2, \mathbf t')$ is straightforward. But we would like to break it up to a block matrix
\begin{equation}
\left(  \begin{array}{cc} \mathcal A_{11} & \mathcal A_{12}\\
\mathcal A_{21}  & \mathcal A_{22} \end{array}\right)
\end{equation}
 where $\mathcal A_{ij}$ are the following Jacobian matrices:

(a)\begin{equation} \mathcal A_{11}={\partial \biggl(f_3(c(t_3)), f_3(c(t_4)), \cdots, f_3(c(t_{5d}))\biggr)\over
\partial (\theta_0^2, \cdots, \hat \theta_1^1, \cdots, \theta_2^d, \epsilon_1, \cdots, \epsilon_{2d})}
\end{equation}

(b) \begin{equation}\begin{array}{c}  \mathcal A_{12}\\
\|\\
{\partial \biggl(f_3(c(t_{5d+1})), f_2(c(t_1)), f_2(c(t_2)), f_1(c(t_1)), f_1(c(t_2)), f_0(c(t_1)), f_0(c(t_2))\biggr)\over
\partial (\theta_0^1, \theta_1^1, r_0, r_1, r_2, r_3, \xi)}\end{array}
\end{equation}

(c) \begin{equation}\begin{array}{c}  \mathcal A_{21}\\
\|\\
{\partial \biggl(f_3(c(t_{5d+1})), f_2(c(t_1)), f_2(c(t_2)), f_1(c(t_1)), f_1(c(t_2)), f_0(c(t_1)), f_0(c(t_2))\biggr))\over
\partial (\theta_0^2, \cdots, \hat \theta_1^1, \cdots, \theta_2^d, \epsilon_1, \cdots, \epsilon_{2d})}\end{array}
\end{equation}

(d) \begin{equation}\begin{array}{c}  \mathcal A_{22}\\
\|\\
{\partial \biggl(f_3(c(t_{5d+1})), f_2(c(t_1)), f_2(c(t_2)), f_1(c(t_1)), f_1(c(t_2)), f_0(c(t_1)), f_0(c(t_2))\biggr)\over
\partial (\theta_0^1, \theta_1^1, r_0, r_1, r_2, r_3, \xi))}\end{array}
\end{equation}

Using lemma 3.3, part (a), $\mathcal A_{11}|_{c_g}$ is a non-zero diagonal matrix and $$\mathcal A_{12}|_{c_g}=0.$$
Therefore it suffices to show 
\begin{equation}
det(\mathcal A_{22})\neq 0.
\end{equation}

Next we apply the ``break-up trick", i.e. we'll change the parameters that determine the matrix, but the change will not 
effect its non-degeneracy \footnote{ This is a trick because the change can't be made before the matrix is broken or reduced. 
For example the change of coordinate's system we are going to make below can't be applied to the entire matrix
$\mathcal A(C_M, f_0, f_1, f_2, \mathbf t')$.}.    
First we change the coordinates to $C_M''$. More precisely we consider
\begin{equation}\begin{array}{c}
\mathcal A_{22}(C_M'')\\
\|\\
{\partial \biggl(f_3(c(t_{5d+1})), f_2(c(t_1)), f_2(c(t_2)), f_1(c(t_1)), f_1(c(t_2)), f_0(c(t_1)), f_0(c(t_2))\biggr)\over
\partial (\theta_0^1, \theta_1^1, r_0, r_1, r_2, r_3, r_4))}\end{array}
\end{equation}

Denote the original $\mathcal A_{22}$ by $\mathcal A_{22}(C_M')$. 
Because at $c_g$, the transformation of coordinates
gives a relation
\begin{equation}
\mathcal A_{22}(C_M')=\mathcal A_{22}(C_M') D_7,
\end{equation}
where $D_7$ is a non-degenerate $7\times 7$ triangular matrix, 
 it suffices to show  $$\mathcal A_{22}(C_M'')$$
is non-degenerate. Notice $t_{5d+1}$ is generic on $\mathbf P^1$. The genericity of $q$ makes
curve in $\mathbb C^7$,
\begin{equation}
({\partial f_3(c(t_{5d+1}))\over \partial \theta_0^1}, {\partial f_3(c(t_{5d+1}))\over \partial \theta_1^1}, {\partial f_3(c(t_{5d+1}))\over \partial r_0}, \cdots, {\partial f_3(c(t_{5d+1}))\over \partial r_4})
\end{equation}
 span the entire space $\mathbb C^7$. This means the first row vector of $$\mathcal A_{22}(C_M'')$$
is generic with respect to other 6 row vectors. 
Hence it suffices for us to show the Jacobian matrix 
\begin{equation}\begin{array}{c}
\mathcal B(c_g)\\
\|\\
{\partial(f_2(c(t_1)), f_2(c(t_2)), f_1(c(t_1)), f_1(c(t_2)), f_0(c(t_1)), f_0(c(t_2)))\over \partial (
(\theta_0^1, \theta_1^1, r_1, r_2, r_3, r_4))}\end{array}
\end{equation} is non degenerate (the column of partial derivatives with respect to $r_0$ is eliminated).
 Now we use  the ``break-up trick" again. This time we change the point $c_g$. To show $\mathcal B(c_g)$ is non degenerate for a generic $c_g\in P(\Gamma_{\mathbb L})$, it suffices to show it is
non-degenerate for a special $c_g\in P(\Gamma_{\mathbb L})$.  To do that, we let $\mathbb L_1$  be an  open set of pencil through $f_0, f_2$.
 $P(\Gamma_{\mathbb L_1})$ be as defined in lemmas 3.2, 3.3.  We choose $P(\Gamma_{\mathbb L_1})$ to be irreducible, and to 
be contained in $P(\Gamma_{\mathbb L})$ for generic  $q$ (simultaneously). 
So a generic point $c_g^1=(c_1^0, \cdots, c_1^4)\in  P(\Gamma_{\mathbb L_1})$ lies in $M_d$.
Use the same notations for $C_M', C_M''$ values of $c_g^1$  as in definition 3.2. 
Because $q$ is generic with respect to 1st, 2nd, 5th and 6th rows, two middle rows 
 \begin{equation}\begin{array}{c}
\biggl({\partial f_1(c(t_1))\over \partial \theta_0^1}, {\partial f_1(c(t_1))\over \partial \theta_1^1}, {\partial f_1(c(t_1))\over \partial r_1}, \cdots,   {\partial f_1(c(t_1))\over \partial r_4} \biggr)\\
\biggl(({\partial f_1(c(t_2))\over \partial \theta_0^1}, {\partial f_1(c(t_2))\over \partial \theta_1^1}, {\partial f_1(c(t_2))\over \partial r_1}, \cdots,   {\partial f_1(c(t_2))\over \partial r_4} \biggr)
\end{array}\end{equation}
must be generic in $\mathbb C^6$ with respect to 1st, 2nd, 5th and 6th rows. Then we reduce $\mathcal B(c_g)$  to show
\begin{equation}
Jac(f_0, c_g^1)={\partial (f_3(c(t_1)), f_3(c(t_1)), f_0(c(t_1)), f_0(c(t_2)))\over \partial (\theta_0^1, r_2, r_3, r_4)}
\end{equation}
is non-degenerate.  Finally we write down the matrix $Jac(f_0, c_g^1)$,
\begin{equation} \begin{array}{c} Jac(f_0, c_g^1)\\
\|\\
\lambda \left(\begin{array}{cccccc}
{1\over t_1-\tilde\theta_0^1} &1  & 1 & 1  \\
{1\over t_2-\tilde\theta_0^1} &1 & 1& 1 \\
 {\partial  f_0(c_g^1 (t_1))\over \partial \theta_0^1}  & (z_2{\partial f_0 \over \partial z_2})|_{c_g^1(t_1)} &
 (z_{3}{\partial f_0 \over \partial z_{3}})|_{c_g^1(t_1)}& (z_4{\partial f_0 \over \partial z_4})|_{c_g^1(t_1)}\\
{\partial f_0(c_g^1 (t_2))\over \partial \theta_0^1} & (z_2{\partial f_0 \over \partial z_2})|_{c_g^1(t_2)} &
 (z_{3}{\partial f_0 \over \partial z_{3}})|_{c_g^1(t_2)}& (z_4{\partial f_0 \over \partial z_4})|_{c_g^1(t_2)}
\end{array}\right), \end{array}\end{equation}
where $\lambda $ is a non-zero complex number.
We further compute to have
\begin{equation}\begin{array}{c}  Jac(f_0, c_g^1)\\
\|\\
\lambda ({1\over t_1-\tilde\theta_0^1}-{1\over t_2-\tilde\theta_0^1} )\left (\begin{array}{ccc}
1&1&1\\  (z_2{\partial  f_0\over \partial z_2})|_{c_g^1(t_1)} &
(z_3{\partial  f_0\over \partial z_3})|_{c_g^1(t_1)} &(z_4{\partial  f_0\over \partial z_4})|_{c_g^1(t_1)} \\
(z_2{\partial  f_0\over \partial z_2})|_{c_g^1(t_2)} &
(z_3{\partial  f_0\over \partial z_3})|_{c_g^1(t_2)} &(z_4{\partial  f_0\over \partial z_4})|_{c_g^1(t_2)} 
\end{array}\right). \end{array}
\end{equation}

Since $t_1, t_2$ are only required to satisfy one equation (5.2), by the genericity of $q$, we may assume $(t_1, t_2)\in \mathbb C^2$ is
generic.  It suffices to prove that 

$$  
\left|\begin{array}{ccc}
1&1&1\\  (z_2{\partial  f_0\over \partial z_2})|_{c_2(t_1)} &
(z_3{\partial  f_0\over \partial z_3})|_{c_2(t_1)} &(z_4{\partial  f_0\over \partial z_4})|_{c_2(t_1)} \\
(z_2{\partial  f_0\over \partial z_2})|_{c_2(t_2)} &
(z_3{\partial  f_0\over \partial z_3})|_{c_2(t_2)} &(z_4{\partial  f_0\over \partial z_4})|_{c_2(t_2)} 
\end{array}\right|\neq 0 
$$
for any generic $f_0$ and $c_2$ that has no multiple zeros with coordinates planes $\{z_i=0\}$. 

Let $\Sigma$ be an open subvariety
$$\{c\in M: zeros \ of \ c_i(t)=0\ are \ distinct, i=1, \cdots, 4\}.$$
Consider the family of rational maps

$$V_{f}=\{ c\in \Sigma: Jac(f, c)=0 \}.$$
Notice by the definition $V_{f}$ is a subvariety of $\Sigma$.
Next we consider the fibre $V_{f_{Fe}}$ where $f_{Fe}$ is the Fermat quintic
$$f_{Fe}=z_0^5+\cdots+z_4^5.$$
It is obvious $V_{f_{Fe}}$ is empty. Hence $V_{f}$ is empty for generic $f$.  This shows that 
\begin{equation} \left|\begin{array}{ccc}
1&1&1\\  (z_2{\partial  f_0)\over \partial z_2})|_{c_2(t_1)} &
(z_3{\partial  f_0)\over \partial z_3})|_{c_2(t_1)} &(z_4{\partial  f_0)\over \partial z_4})|_{c_2(t_1)} \\
(z_2{\partial  f_0)\over \partial z_2})|_{c_2(t_2)} &
(z_3{\partial  f_0)\over \partial z_3})|_{c_2(t_2)} &(z_4{\partial  f_0)\over \partial z_4})|_{c_2(t_2)} 
\end{array}\right|\neq 0. 
\end{equation}
Therefore 
$$Jac(f_0, c_2)\neq 0.$$

We complete the proof of lemma 4.2. 

\end{proof}

\subsection{Ranks of differential sheaves}
\quad\smallskip

\bigskip

\begin{proof} of proposition 1.4: 
Let $\mathcal N$ be the submodule of global sections, $H^0(\Omega_{M})$
generated by
elements
\begin{equation} \phi_i= d
\left|  \begin{array}{ccc} f_2(c(t_1)) & f_1(c(t_1)) & f_0(c(t_1))\\
f_2(c(t_2)) & f_1(c(t_2)) & f_0(c(t_2))\\
f_2(c(t_i)) & f_1(c(t_i)) & f_0(c(t_i))\end{array}\right|\end{equation}
for $i=3, \cdots, 5d+1$. 
Recall that 
$$\left|  \begin{array}{ccc} f_2(c(t_1)) & f_1(c(t_1)) & f_0(c(t_1))\\
f_2(c(t_2)) & f_1(c(t_2)) & f_0(c(t_2))\\
f_2(c(t_i)) & f_1(c(t_i)) & f_0(c(t_i))\end{array}\right|=0,$$
for $i=3, \cdots, 5d+1$
define the scheme $P(\Gamma_{\mathbb L})$ for a small $\mathbb L$.  By   proposition 8.12 in [7], II, \begin{equation}
\widetilde {({H^0(\Omega_{M})\over \mathcal N})}\otimes \mathcal O_{P(\Gamma_{\mathbb L})}
\simeq \Omega_{P(\Gamma_{\mathbb L})},
\end{equation}
where $\widetilde {(\cdot)}$ denotes the sheaf associated to the module $(\cdot)$.

Therefore

\begin{equation}\begin{array}{cc} &
 {({H^0(\Omega_{M})\otimes k(c_g)\over \mathcal N\otimes k(c_g)})}
\simeq \Omega_{P(\Gamma_{\mathbb L})}\otimes k(c_g)\\&
=(\Omega_{P(\Gamma_{\mathbb L})})|_{(\{c_g\})},\end{array}\end{equation}
where $k(c_g)=\mathbb C$ is the residue field at generic $$c_g\in P(\Gamma_{\mathbb L}).$$ 
Notice two sides of (4.31) are finitely dimensional  linear spaces over $\mathbb C$.
\begin{equation}\begin{array}{cc} &
dim_{\mathbb C} ((\Omega_{P(\Gamma_{\mathbb L})})|_{(\{c_g\})})\\
&=dim_{\mathbb C}(H^0(\Omega_{M})\otimes k(c_g))
- dim(\mathcal N\otimes k(c_g))\end{array}
\end{equation}

Since \begin{equation}
dim_{\mathbb C} ((\Omega_{P(\Gamma_{\mathbb L})})|_{(\{c_g\})}))=dim(T_{c_g}P(\Gamma_{\mathbb L}))
\end{equation}

\begin{equation}
dim(T_{c_g}P(\Gamma_{\mathbb L}))=dim(M) -dim(\mathcal N\otimes k(c_g)).
\end{equation} 

By lemma 4.1, 
$$dim(\mathcal N\otimes k(c_g))=deg(\omega(M, \mathbf t)).$$

The proposition 1.4 is proved.

\end{proof}

\bigskip

 \begin{proof} of theorem 1.1.  This is the statement of proposition 1.5. See section 3, [9] for the details.
\end{proof}

\section{Examples\\
--Vainsencher's and Chen's rational curves}

{\bf Example 5.1} (Vainsencher's rational curves)
 \par
This example provides an evidence to theorem 1.1.
In [8], Vainsencher constructed irreducible, degree 5, nodal curves $C_0$ on a generic quintic $f_0$ by taking plane sections of the quintic. 
Let $c_0$ be its normalization. By
our theorem 1.1, $c_0$ is an immersion and 
\begin{equation}
N_{c_0/X_0}\simeq \mathcal O_{\mathbf P^1}(-1)\oplus \mathcal O_{\mathbf P^1}(-1).
\end{equation}

Indeed these were proved by Cox and Katz in [4], by using a different method. Their method is based on Clemens' deformation idea.   Their understanding of $c_0$ on $f_0$  was achieved by a concrete construction 
of special $c_0$, $f_0$ and by using a computer program for the
last verification of the $26\times 30$ matrix.  It is easy to check that the rational maps $c_0$ they constructed are immersions.  \par
Furthermore
our result shows 
\begin{equation} dim(T_{c_0}\Gamma_{f_0})=4.
\end{equation}
Because of the equation (5.1), $C_0$ can't deform in $f_0$. Thus
$\Gamma_{f_0}$ consists of multiple orbits isomorphic to $GL(2)(c_0)$.   Theorem 1.1 also shows that there will not be any scheme-theoretical multiplicity associated to the orbits. However the number of these orbits is not accessible because the degree of each orbit in $\mathbf P(M)$
could be different. This number is related to  Gromov-Witten invariants.

\bigskip

{\bf Example 5.2} (Chen's rational curves) \par

This is an example on $K$-3 surfaces. 
In [1], Chen constructed  nodal rational curves $C_0$ of degree $4d$ for each natural number $d$,  that lie on the 
generic hypersurfaces $f_0$ of degree $4$ in $\mathbf P^3$ ($f_0$ is a K-3 surface).   At first we may have an impression that this is against our intuition. 
Because it is similar to rational curves on generic  quintic threefolds that we can have naive counting: on a generic quartic  hypersurface $f_0$ of $\mathbf P^3$ , there will be $4d+1$ conditions imposed the rational curves on $f_0$, while the dimension of the moduli space of rational curves in $\mathbf P^3$ (modulo $PGL(2)$ action) is only $4d$. Thus the naive counting concludes that there will not be any rational curves on $f_0$.  But it was proved by Mori, Mukai, etc.,  and  Chen ([1]) that rational curves on $f_0$ exist and they are all nodal. Our proof is closely related to this counting, and our construction of $\omega(M, \mathbf t)$ can be carried out in $\mathbf P^3$ for Chen's case.  
But theorem 1.1 does not hold because proposition 1.3 fails.  This failure is not expected by the naive dimension count, but it is a reminder of a fact that  the  generic quartics are not generic in the moduli space of complex  structures.    \par
 
Chen's  construction has a similar
flavor of Vainsencher's rational curves  above. They were obtained by taking hyperplane sections of $K$-3 surfaces. 
Intrinsically Vainsecher's and Chen's rational curves look similar.  For instance they are all plane sections, and are all immersed,  nodal rational curves.   So what invariant distinguishes one from the other?  Section 4 shows that this invariant may not be the invariant of the intrinsic rational curves, it  addresses the structure of the moduli space of rational curves for underlined families of varieties. 
 More specifically, it is deduced from the differential form $\omega(M, \mathbf t)$ (defined in (1.10) ). The $\omega$ itself is not a moduli invariant, but the 
zero locus $\{\omega(M, \mathbf t)=0\}$ is, and furthermore $\{\omega(M, \mathbf t)=0\}$ is independent of generic $t_i, i=1, \cdots, 5d+1$. In Chen's situation, $\omega(M, \mathbf t)$ turns out to be  identically zero on  $P(\Gamma_{\mathbb L})$,  but in Vainsencher's  it is not.   Beyond Chen's cases, it is not clear that which homology classes of rational curves would have or would not have vanishing $\omega(M, \mathbf t)$.

\end{document}